\documentclass[reqno,a4paper,12pt]{amsart}

\usepackage[full]{textcomp}
\usepackage{fourier}
\usepackage{latexsym}

\usepackage[pdftex]{graphicx}
\graphicspath{{IMAGES/}{IMAGES1/}{./}}

\usepackage[backref=page]{hyperref} 

 \usepackage{amscd}
\usepackage{booktabs}
\newcommand{\ra}[1]{\renewcommand{\arraystretch}{#1}}

\newcommand{\ww}{{\mathfrak w}}
\newcommand{\vv}{{\mathfrak v}}
\newcommand{\dd}{{\mathfrak d}}
\newcommand{\hh}{{\mathfrak h}}
\def\FF{\mathbb{F}}
\def\NN{\mathbb{N}}
\def\FF{\mathbb{F}}
\def\ZZ{\mathbb{Z}}

\def\PP{\mathcal P}

\newtheorem{theorem}{Theorem}
\newtheorem{proposition}{Proposition}
\newtheorem{corollary}{Corollary}

\theoremstyle{remark}

\begin{document}
%


\title{A game with divisors and absolute differences of exponents} 
\author{Cristian Cobeli and Alexandru Zaharescu}

\address{
\indent CC \textit{and} AZ: 
\emph{Simion Stoilow},
 Institute of Mathematics of the Romanian Academy,  \newline
\indent
P.O. Box 1-764, 
Bucharest, 70700, Romania.}
\email{cristian.cobeli@imar.ro}

\address{
\indent AZ:
Department of Mathematics, 
University of Illinois at Urbana - Champaign,  \newline
\indent 1409 West Green Street, 
Urbana, IL 61801, USA.
}
\email{zaharesc@illinois.edu}

\keywords{\em  Ducci game, Gilbreath's conjecture, primes game, absolute differences,
Sierpinski triangle.}
\subjclass[2010]{
Primary 11B85; 
Secondary 39A10.
}

\begin{abstract}
In this work, we discuss a number game that develops in a manner similar to
that on which Gilbreath's conjecture on iterated absolute differences between consecutive primes is
formulated.
In our case the action occurs at the exponent level and there, the evolution is reminiscent of
that in a final Ducci game.
We present features of the whole field of the game created by the successive generations, prove
an analog of Gilbreath's conjecture, and raise some open questions.
\end{abstract}

\maketitle

\section{Introduction}
We study a cellular automaton, whose different facets display characteristics that are both
multiplicative and additive in nature.
Essentially, it is defined by an 'atomic transformation' function, which is given by an arithmetical
operation.
To any integers $a$, $b$ we associate the integer
\medskip

\begin{equation*}
   Z(a, b):=\frac{ab}{\big(\gcd(a, b)\big)^2}\,.
\end{equation*}
\medskip
\smallskip

\noindent
Applying this transformation on any two neighbor elements of a given 
sequence of integer numbers, we obtain a new 
$Z(\cdot, \cdot)$-generated sequence. We write it on the next row, such that each term is placed in
the middle, below its parents, like in a contemporary presentation of the Pascal arithmetic
triangle. One difference, however, lies in the fact that in our cases of interest the first row is
long.

\begin{center}
\begin{figure}[ht]
\vspace*{7mm}
\centering
   \includegraphics[height=0.89\textwidth,width=0.53\textwidth,angle=-90]{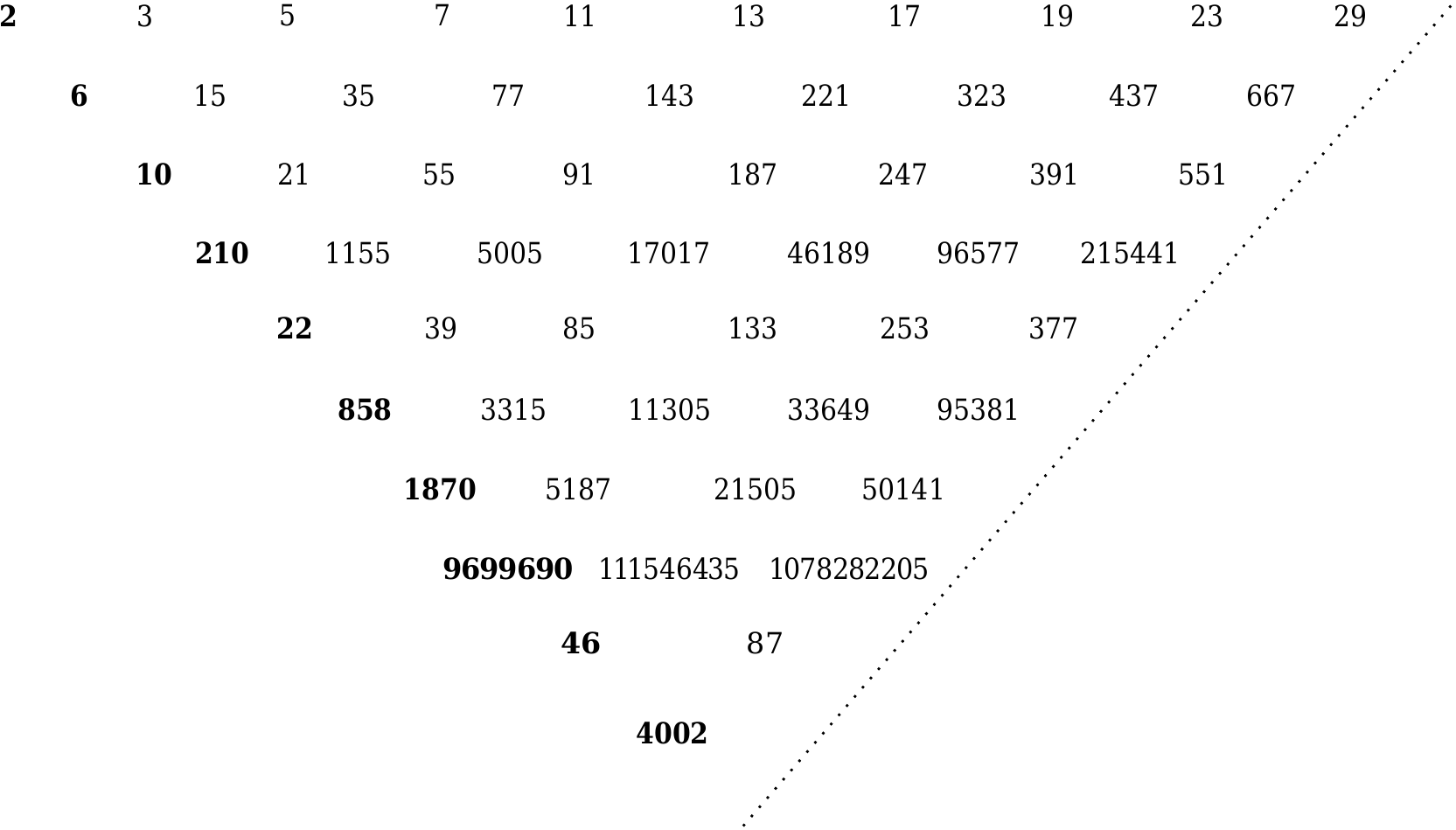}
\caption{The triangular slice of order  $10$ generated by iterated applications of 
$Z(\cdot, \cdot)$ on the sequence of prime numbers.}
 \label{FigureZprimes10}
 \end{figure}
\end{center}

Our subject of study is the construction that starts 
with the sequence of prime numbers on the first row. 
We denote by $a_{m,n}$ the $n$-th element from the $m$-th row of this
table, for $m\ge 0$, $n\ge 1$.  
Looking at the initial part of the honeycomb table of numbers thus generated, we get a triangle such
as the one displayed in Figure~\ref{FigureZprimes10}.
During the involved calculations one finds, for example:

\bigskip

\begin{equation*}
\begin{split}
   Z(46189,\  96577) &= Z(11\cdot 13\cdot 17\cdot 19, \ \; 13\cdot 17\cdot 19\cdot 23)\\
&=\frac{11\cdot 13^2\cdot 17^2\cdot 19^2\cdot 23}{13^2\cdot 17^2\cdot 19^2}
=11\cdot 23
=253.
\end{split}
\end{equation*}

{
\newcommand{\mc}[3]{\multicolumn{#1}{#2}{#3}}
\begin{center}
\begin{table}[ht]
\ra{1.1}
\caption{The first $16$ entries on the left edge and their prime factorization.}
\medskip 

\centering
\begin{tabular}{@{}cll@{}}
\toprule
$m$ & $d_m$ & factorization of $d_m$\\
\midrule
$0$ & $2$ & $2$\\
$1$ & $6$  & $2\cdot 3$\\
$2$ & $10$  & $2\cdot 5$\\
$3$ & $210$  & $2\cdot 3\cdot 5\cdot 7$\\
$4$ & $22$  & $2\cdot 11$\\
$5$ & $858$  & $2\cdot 3\cdot 11\cdot 13$\\
$6$ & $1870$  & $2\cdot 5\cdot 11\cdot 17$\\
$7$ & $9699690$  & $2\cdot 3\cdot 5\cdot 7\cdot 11\cdot 13\cdot 17\cdot 19$\\
$8$ & $46$ & $2\cdot 23$\\
$9$ & $4002$  & $2\cdot 3\cdot 23\cdot 29$\\
$10$ & $7130$ & $2\cdot 5\cdot 23\cdot 31$\\
$11$ & $160660290$  & $2\cdot 3\cdot 5\cdot 7\cdot 23\cdot 29\cdot 31\cdot 37$\\
$12$ & $20746$  & $2\cdot 11\cdot 23\cdot 41$\\
$13$ & $1008940218$  & $2\cdot 3\cdot 11\cdot 13\cdot 23\cdot 29\cdot 41\cdot 43$\\
$14$ & $2569288370$  & $2\cdot 5\cdot 11\cdot 17\cdot 23\cdot 31\cdot 41\cdot 47$\\
$15$ & $32589158477190044730$  & $2\cdot 3\cdot 5 \cdot 7\cdot 11\cdot 13\cdot 17\cdot 19\cdot
23\cdot 29\cdot 31\cdot 37\cdot 41\cdot 43\cdot 47\cdot 53$\\
\bottomrule
\end{tabular}\label{Table1}
\end{table}
\end{center}
}%

A construction similar to that in Figure~\ref{FigureZprimes10}, but with $Z(\cdot,\cdot)$
replaced by the absolute difference function leads to
a table whose appearance seems simpler (see Guy~\cite[Problem A10]{Guy1994}). 
In that case the numbers quickly get smaller and one notices that
the numbers on the left edge of the triangle are equal to one. Proth~\cite{Pro1878} thought
that he proved that fact, but his argument was incorrect and the problem remained open. Later
Gilbreath rediscovered the property
and the statement is known as Glibreath's Conjecture. 
Odlyzko~\cite{Odl1993} verified the conjecture for triangles of size $<\pi(10^{13})$.
The apparent simplicity of the rows of higher generations hides deep problems involving gaps of
different orders between primes, on which 
current results 
only
scrape the surface of a field that remains untouchable by known techniques.
Although, Odlyzko~\cite{Odl1993} observed that there are chances to overcome the
involved difficulties. Reasoning
heuristically, he concluded that the taming property of the numbers on the left edge of such a
construction should be generally valid. Indeed, he showed that due to iterativity, 
this should be a characteristic property of lots of other 
sequences, provided their gaps to neighbor elements are random and relatively small.

On the contrary, in our case,
one sees that the numbers become large, growing in augmented
waves. Notice that by definition it follows that they are all square free.
Looking for an analogue of Gilbreath's conjecture, we consider 
the prime factorization of the
first element of each row of the construction whose first part is
presented in Figure~\ref{FigureZprimes10}. 
Denote by $d_m$ the first element of row $m$, for $m\ge 0$. The factorization into primes of the
first
few of them is
given in Table~\ref{Table1}. Analyzing this table, we remark a striking 
phenomenon: the number of prime factors of each of these numbers is a power of  $2$.
The following theorem that will be proved twice, in Sections~\ref{SectionTh1} and
\ref{SectionTh1Th2}, establishes this fact.

\begin{theorem}\label{Theorem1}
   The number of prime factors of $d_m$ is a power of~\ $2$, for any $m\ge 0$.
\end{theorem}

Here the sequence of exponents is itself remarkable.
For any $m\ge 0$, let 
\begin{equation*}
     \delta_m:=\frac{\log \big(\text{the number of prime factors of $d_m$}\big)}{\log 2}\,.
\end{equation*}
Sequence $\{\delta_m\}_{m\ge 0}$ is special and has interesting properties.
For example, it is a proper subsequence of itself, and as such
it is a \emph{fractal sequence} (cf. Kimberling~\cite{Kim1995}). 
Indeed, the subsequence of terms having even rank coincides with the
original sequence. And this fact occurs for sub-subsequences,  also, again and again, ad infinitum.
Precisely this follows by 
\begin{equation}\label{deltasymmetry}
\delta_{t2^s}=\delta_t, \quad\text{for any $0\le t $ and $1\le s$.}   
\end{equation}
We do not know a closed formula for $\delta_m$, but
we know that it verifies many other recursive
formulas. Thus, given an initial group of $2^s$ terms for any $s\ge 1$, we obtain the next group
of $2^s$ terms just by incrementing by $1$ all the known terms, that is, 
\begin{equation*}
\delta_{2^s+j}=\delta_j+1, \quad\text{for any $0\le j \le 2^s-1$ and $1\le s$.}  
\end{equation*}
\begin{center}
\begin{figure}[ht]
\centering
\hfill
   \includegraphics[width=0.48\textwidth]{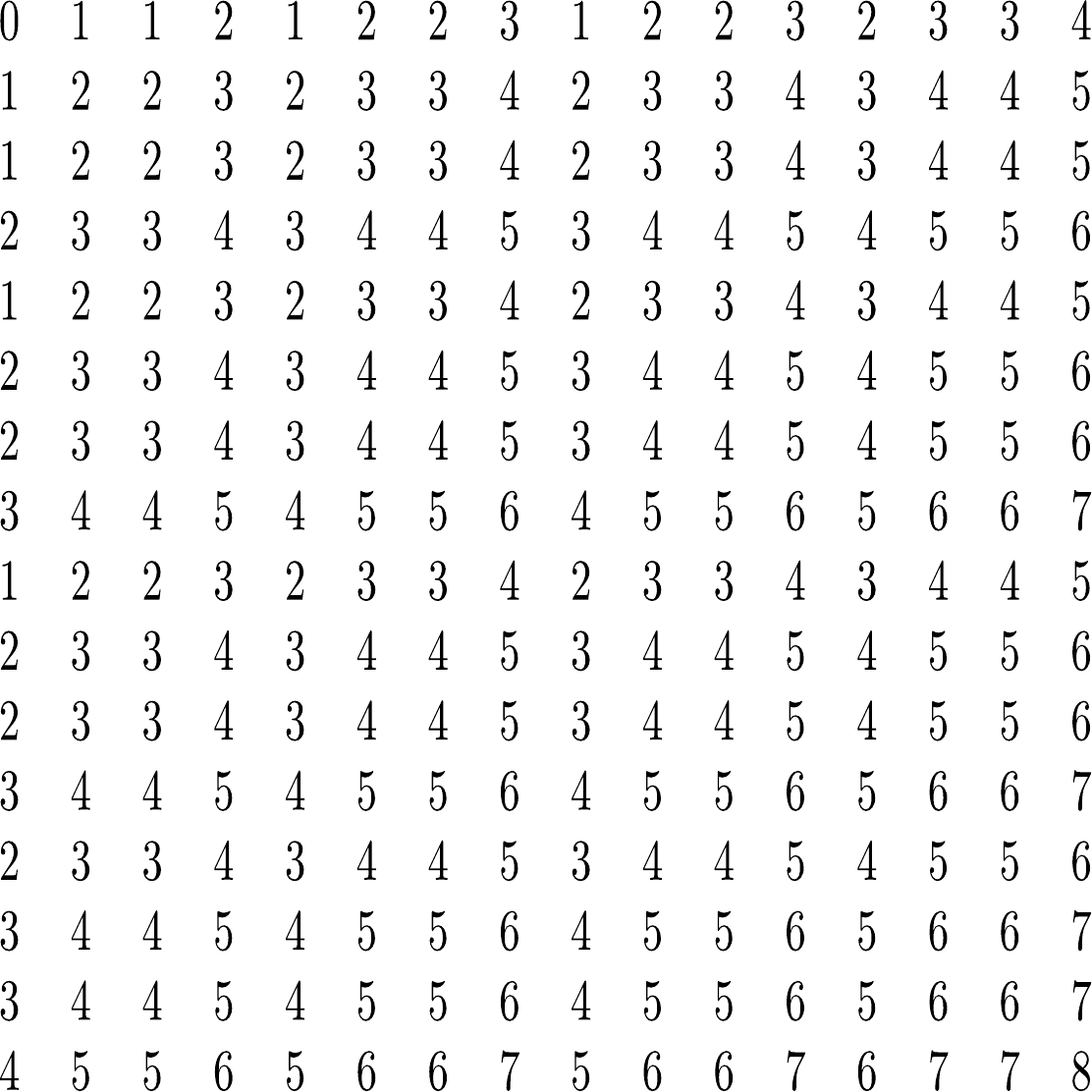}
\hfill\hfill
   \includegraphics[width=0.48\textwidth]{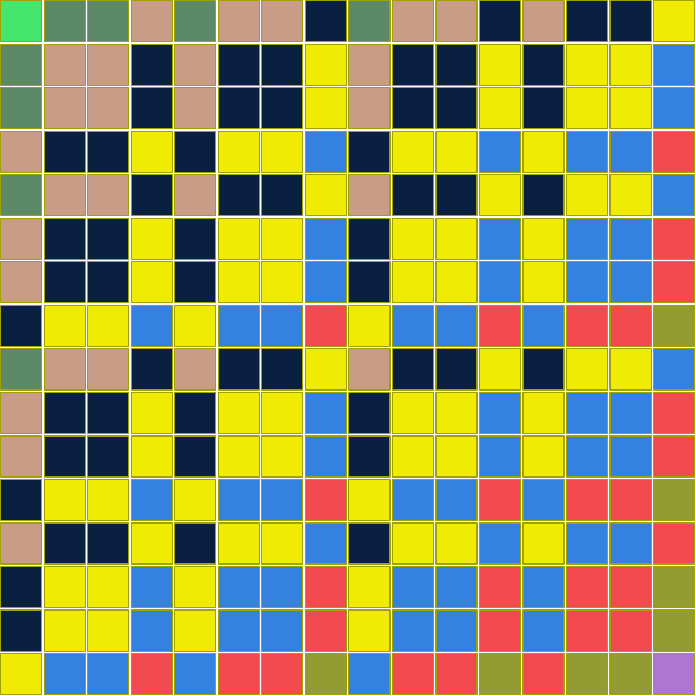}
\hfill\mbox{}
\caption{Sequence of exponents $\{\delta_m\}_{m\ge 0}$ presented in rows of length $16$: as
numbers on the left side and as suitable square cells, filled with corresponding colors, which were
chosen randomly.
}
 \label{Figure16}
 \end{figure}
\end{center}

Computer experiments suggest that the sequence $\{\delta_m\}_{m\ge 0}$ has lots of
$2$-adic symmetries, both local and global, but their proofs are not the object of this work.
However, here we emphasize some of these properties that are related to our automaton, presenting
the sequence in different forms. 
On the left side of Figure~\ref{Figure16} we wrote the first part of the sequence with $16$ elements
on each row. Thus $\delta_0=0$, $\delta_1=1$, \ldots, $\delta_{15}=4$. Next on the second row, we
put: $\delta_{16}=1$, $\delta_{17}=2$, \ldots, $\delta_{31}=5$, and so on. To get a square table,
only the first $16$ rows are shown, but the sequence continues likewise.
A visual representation of the matrix is given on the right side of Figure~\ref{Figure16}. There,
nine distinct random colors are chosen for each of the possible values of $\delta_m$. 
We remark that this square is symmetric with respect to the north-west--south-east diagonal. 
Observe that relation \eqref{deltasymmetry} is merely a particular instance of this property.
%
\begin{center}
\begin{figure}[ht]
\centering
\hfill
\begin{minipage}{0.48\textwidth}
\centering
        \includegraphics[width=1\textwidth]{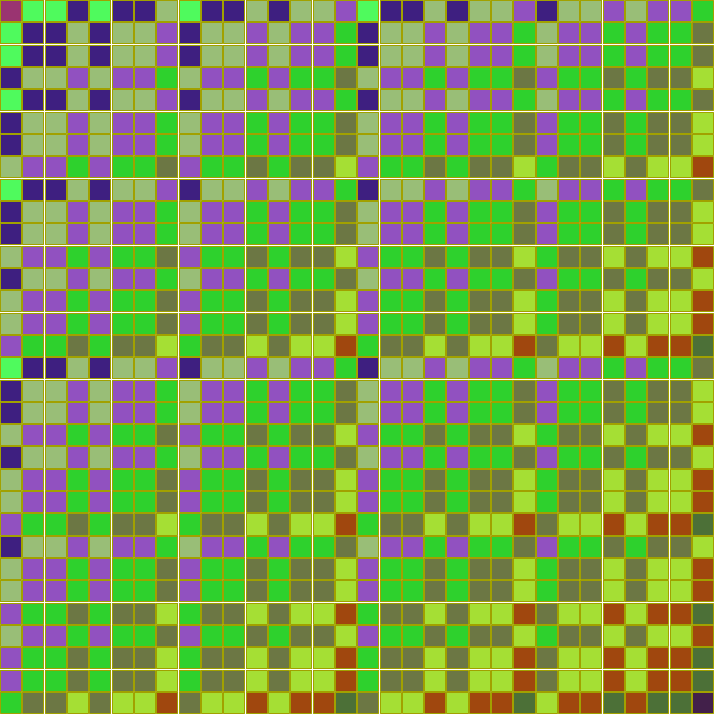}%
\end{minipage}
\hfill\hfill
\begin{minipage}{0.48\textwidth}
\centering
        \includegraphics[width=1\textwidth, angle=-90]{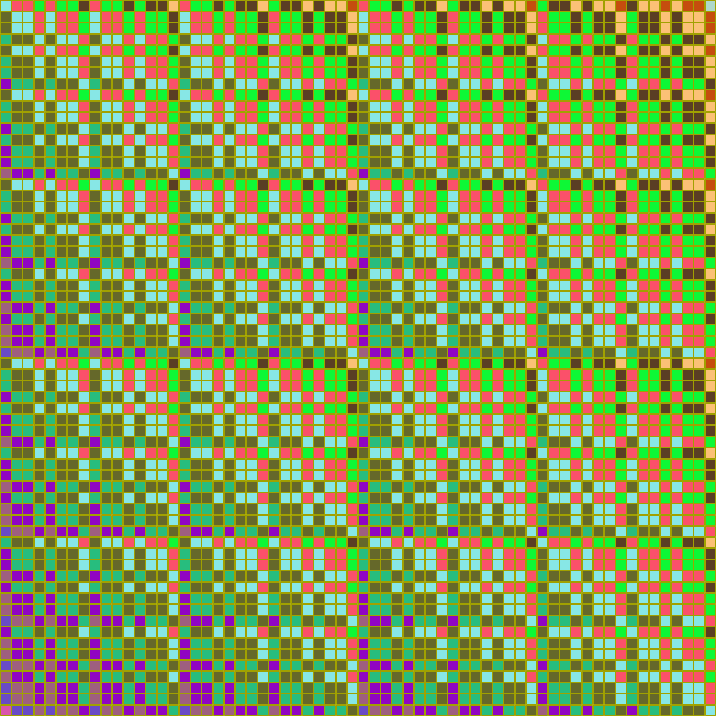}%
\end{minipage}
\hfill\mbox{}
\caption{Representations of the sequence $\{\delta_m\}_{m\ge 0}$ in tables  having $32$ and $64$
cell in each row, respectively. In each case, the colors that corresponds to the numbers 
$\delta_m$ are chosen  randomly.
}
 \label{Figure32-64}
 \end{figure}
\end{center}
%
Another 'complementary symmetry' is featured by the table
 with respect to the south\-west--north-east diagonal.  We observe
 that all the elements on this diagonal have the same color, that of $\delta_{15}=4$.
We also see that, if $\delta_{r}$
and $\delta_{s}$ are on  squares symmetric  with respect to the south-west--north-east diagonal,
then
$\delta_{r} +\delta_{s}=8$. 

Furthermore, the square from Figure~\ref{Figure16} may also be generated starting with the basic 
\textit{'number motif'} (or \textit{'colored motif'}) of the $4\times 4$ square situated on the
north-west
corner. This motif is repeated with a period of length $4$ again and again along the horizontal
or vertical lines, and also along the $45^\circ$ diagonals. 
Along the way, the only change is the simultaneous  increment or decrement of all the
elements in the number motifs by a specific value.
The values that change the  number motifs evolve $2$-adically
 in a remarkable way. 
To better emphasize these properties, in  Figure~\ref{Figure32-64} we show the analogue squares
generated in a similar way, by writing the sequence $\{\delta_m\}_{m\ge 0}$ in rows of $32$ and 
$64$ at a time, respectively.
Notice also that one may consider as the basic number motif a $2\times 2$ square and also any
larger square of size $2^t\times 2^t$.
These symmetries occur also in analogue packages of the sequence $\{\delta_m\}_{m\ge 0}$ in boxes
embedded in higher dimensional spaces. Despite of this and the simplicity of their construction,
the pictures are on the verge between regularity and randomness.

Sequence $\{\delta_m\}_{m\ge 0}$ has many other interesting properties.
It is sequence {\tt A000120} from Sloane's On-line Encyclopedia of Integer Sequences and it is
related to other remarkable sequences, cf.~\cite{OEIS} and the references therein.
They are mainly associated with binary representations, because $\delta_m$ equals the number of
$1$'s
from the $m$-th row of Pascal triangle modulo $2$. 
Further, this links them to Fermat numbers,
constructible polygons and Bernoulli numbers. A short survey on these ties is presented
in~\cite[Section 4]{CZ2013}.

Further inspection of the factorization of entries in Figure~\ref{FigureZprimes10} reveals that
$\omega(a)$, the number of prime factors of $a$, is the same for all integers $a$ situated on the
same horizontal line. 
This is the object of our next theorem that will be proved in Section~\ref{SectionTh1Th2}.
\begin{theorem}\label{Theorem2}
   We have $\omega(d_m)=\omega(a_{m,n})$, for all $m\ge 0$, $n\ge 1$.
\end{theorem}

Deeper properties of the 
evolution function $Z(\cdot, \cdot)$
may be investigated in other contexts,
where  $Z(\cdot, \cdot)$ is applied to integers
that are not necessarily square free as in the 
present case. Various interesting questions arise for instance
when the role of the starting sequence (played by the primes in this article) is
taken by the sequence of natural numbers, or by others, such as the Fibonacci sequence.
If the first row is a line of binomial coefficients
$\binom{n}{0}$, $\binom{n}{1},\dots,\binom{n}{n}$, 
 from Pascal's arithmetic triangle, the generated table is a bounded  triangle. Looking at
their $p$-renderings, one sees striking differences in these pictures
when one of  $p$ or $n$
is kept fixed and the other varies, even slightly (cf.~\cite{CZ2013}).

Working only with integers that are products of distinct primes, enables
one to develop techniques that may be useful in the general case.
In the following sections we elaborate such a method that may be useful in all cases and the main
results stated above will follow, too.


\section{A related game}

In this section we apparently depart from the original  problem, experimenting with a setting that is
more general and more particular at the same time. 
For this, we lift the left edge barrier. This extends the field of the game to also allow in the
new generations entries  from the left side. 
Next, we assume here  that the entries are elements of
$\FF_2=\ZZ/2\ZZ$ and the operation  $Z(\cdot, \cdot)$ is replaced by
the difference in
$\FF_2$. 
This is close to the essence of the end phase of an $n$-absolute difference game
cf. the authors'~\cite{CCZ2000}. In the last century, this automaton was discovered independently,
several times, by different authors.
Their attention was mainly attracted by the fact that  the Ducci game enters into a null cycle when
$n$ is a power of  $2$.
They also studied the length of the period for an arbitrary $n$ and
the number of the necessary iterations before the game enters into a cycle.
For these and further problems and generalizations, the reader is referred to
 Ciamberlini and Marengoni~\cite{CM1937},
Meyers \cite{Mey1982},
Thwaites~\cite{Thw1996a},~\cite{Thw1996b}, 
Pompilli \cite{Pom1996},
Cr\^a\c smaru and the authors'~\cite{CCZ2000},
 Andriychenko and Chamberland \cite{AC2000},
Chamberland \cite{Cha2003a},
 Behn et al.~\cite{BKP2005},
Thomas et al. \cite{CT2004}, \cite{CST2005},
Webb \cite{Web1982},
Ehrlich \cite{Ehr1990},
Ludington-Young \cite{LY1990}, \cite{LY1999},
Calkin et al. \cite{CST2005},
 Lidman and Thomas \cite{LT2007},
 Brown and Merzel \cite{BM2007}
Breuer \cite{Bre2007}, \cite{Bre2010},  
 Caragiu et al. \cite{BC2007},  \cite{CZZ2011} and \cite{CZ2013}.

\begin{center}
\begin{figure}[ht]
\centering
   \includegraphics[scale=0.72,angle=-90]{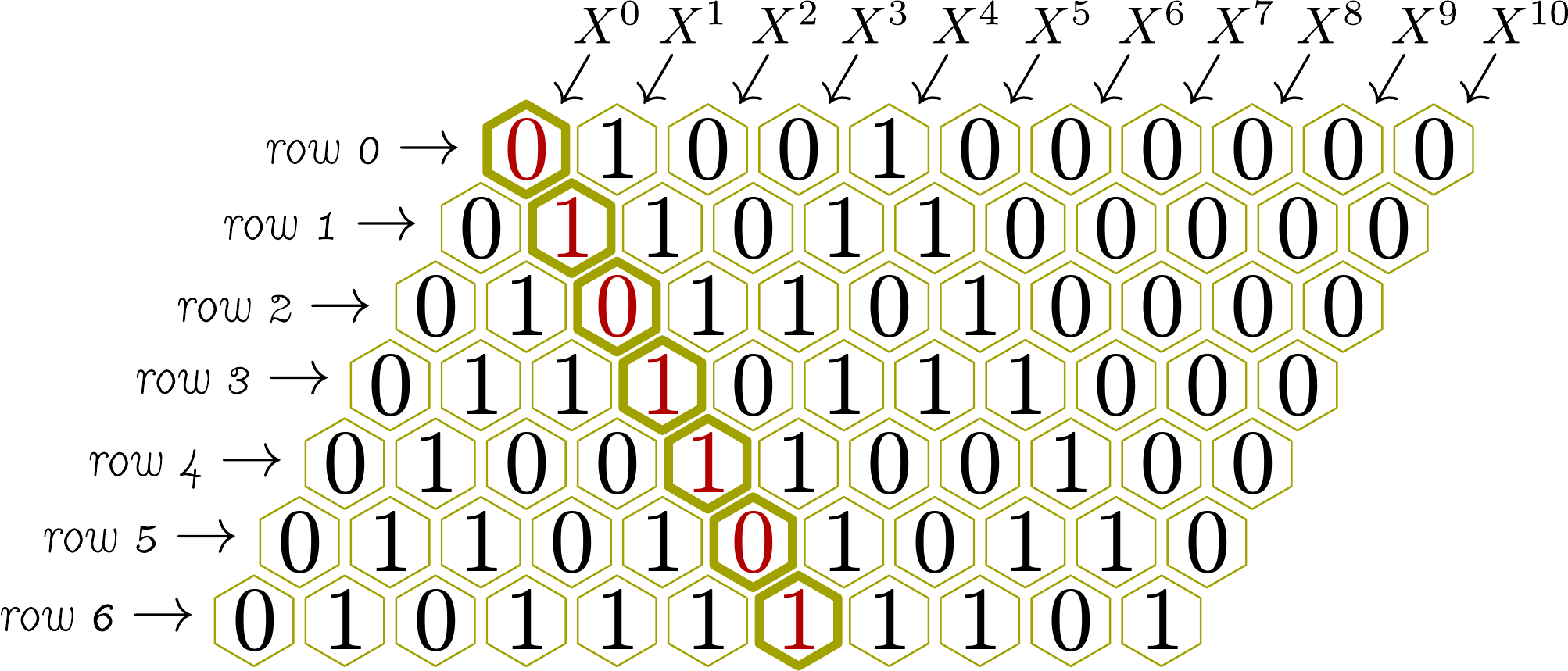}
\caption{The coefficients of the first polynomials of the sequence that starts with $X+X^4$.}
 \label{Figure1}
 \end{figure}
\end{center}

We denote by  $\ww\colon\ \NN \to \FF_2$  the sequence that lists the numbers from the first
line of the game.
Next, we denote by $\psi$ the function that transforms an ascendent sequence $\mathfrak s$
into the succeeding generation $\psi(\mathfrak s)$,
by taking  the 
difference in $\FF_2$ between all neighbor elements. 
For example,  Figure~\ref{Figure1} shows the honeycomb representation of  the generations 
$\psi^{(0)}(\vv)=\vv$, 
$\psi^{(1)}(\vv)$, \ldots,
$\psi^{(6)}(\vv)$, 
where
$\vv(1)=\vv(4)=1$ and $\vv(n)=0$, otherwise. (We denote by $\psi^{(n)}$ the $n$-th composition of
the function $\psi$ by itself.)
There, the highlighted cells are the original positions of the numbers of interest in the
game presented in the introduction and on its left is a new growing triangle.
We denote by $\dd_m=\dd_m(\ww)$, $0\le m$, the $m$-th number situated on that diagonal in the game
that starts
with the sequence $\ww$.
To any sequence $\ww$, we associate the series $P_{\ww}(X)\in\FF_2[[X]]$\footnote{We denote by
$\FF_2[[X]]$ the ring of formal power series, the natural extension of $\FF_2[X]$, the ring of
polynomials  of one variable  over $\FF_2$.} whose coefficients are the
corresponding elements of the sequence $\ww$. If the number of nonzero elements is finite, then
$P_{\ww}(X)$ is a polynomial in $\FF_2[X]$.

For example, the polynomials associated with the first generations shown in Figure~\ref{Figure1}
are:
\begin{equation*}
   \begin{split}
      P_{\vv}(X)&=X+X^4,\\
      P_{\psi(\vv)}(X)&=(X+X^4)(1+X)=X + X^2 + X^4 + X^5,\\
      P_{\psi^{(2)}(\vv)}(X)&=(X+X^4)(1+X)^2=X + X^3 + X^4 + X^6,\\
      P_{\psi^{(3)}(\vv)}(X)&=(X+X^4)(1+X)^3=X + X^2 + X^3 + X^5 + X^6 + X^7,\\
      P_{\psi^{(4)}(\vv)}(X)&=(X+X^4)(1+X)^3=X + X^4 + X^5 + X^8,\\
      P_{\psi^{(5)}(\vv)}(X)&=(X+X^4)(1+X)^3=X + X^2 + X^4 + X^6 + X^8 + X^9,\\
      P_{\psi^{(6)}(\vv)}(X)&=(X+X^4)(1+X)^3=X + X^3 + X^4 + X^5 + X^6 + X^7 + X^8 + X^{10}.\\
   \end{split}
\end{equation*}
The operation of passing from one generation to the next
one corresponds to the
polynomial side to multiplication by $1+X$ in the ring $\FF_2[X]$.
Consequently, finding the $m$-th generation that starts with $\ww$, reduces to finding the
associated series:
\begin{equation*}
\begin{CD}
\ww @>     >> P_{\ww}(X)\\
@VVV			@VV\times(1+X)^m V\\
\psi^{(m)}(\ww)		@>  >> 		P_{\psi^{(m)}(\ww)}(X)
\end{CD}   
\end{equation*}
where we have
\begin{equation}\label{eqpowern}
P_{\psi^{(m)}(\ww)}(X)=P_{\ww}(X)(1+X)^m.
\end{equation}
Relation \eqref{eqpowern} enables one to obtain information on the periodicity of the
generating system.

\begin{center}
\begin{figure}[ht]
\centering
   \includegraphics[scale=0.7,angle=-90]{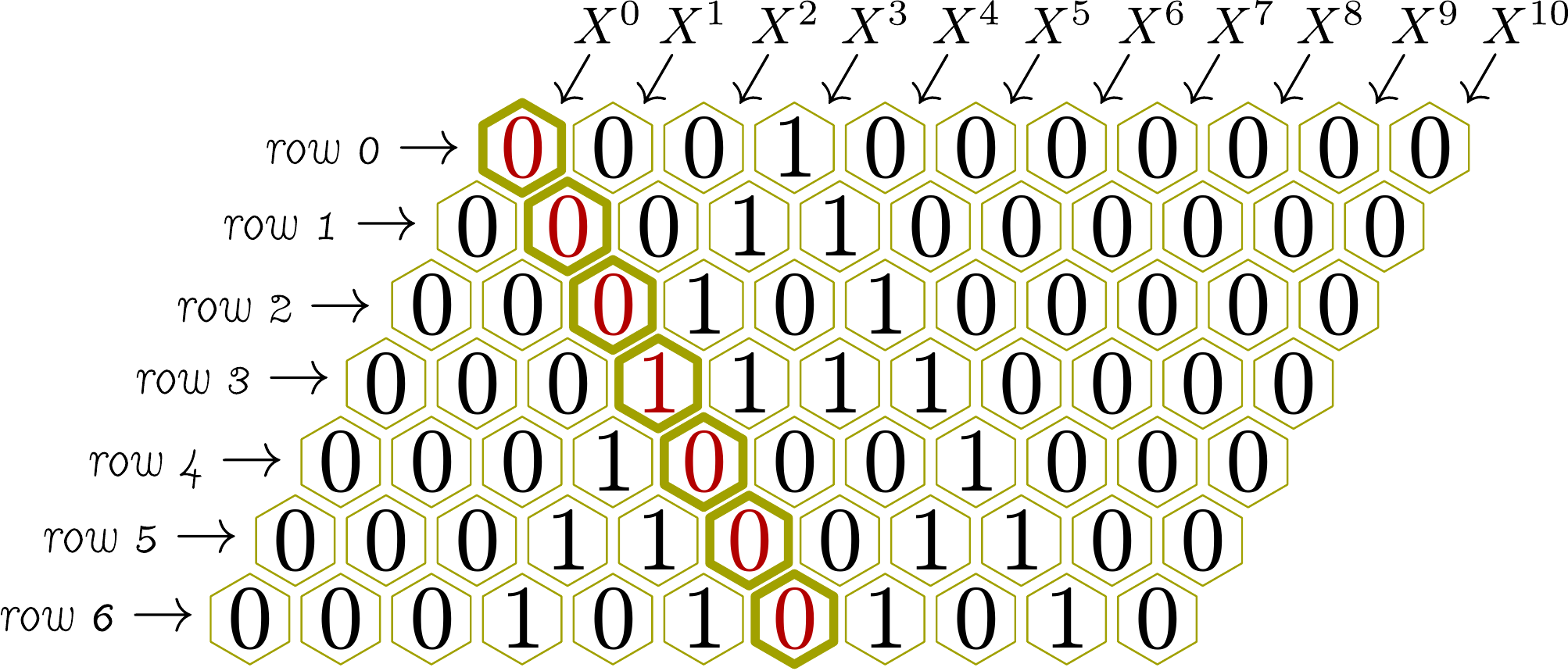}
\caption{The coefficients of the first polynomials of the sequence that starts with $X^3$.}
 \label{Figure2}
 \end{figure}
\end{center}

Denote by $\nu_p(a)$ the exponent of $p$ in the prime factorization of
$a$.
Then, if $0< t <2^s$, $s\ge 1$ and $t_1=t/2^{\nu_2(t)}$, it follows that
$\nu_2\big(2^s-t\big)=\nu_2\big(2^{\nu_2(t)}(2^{s-\nu_2(t)}-t_1)\big)=\nu_2(t)$.
Therefore
\begin{equation*}
   \binom{2^{s}-1}{k}=\prod_{t=1}^{k}\frac{2^s-t}{t} \quad\text{is odd for $0<k < 2^s$,}
\end{equation*}
and $\binom{2^{s}}{k}=0$ in $\FF_2$ for $0< k < 2^s$.
Using this in combination with \eqref{eqpowern}, we obtain the following shifting-enhanced periodicity
feature of the $\FF_2$ game.
%
\begin{proposition}\label{Proposition1}
   For any sequence $\ww$ and any integer $s\ge 1$,  we have
\begin{equation*}\label{eqpower2nB}
P_{\psi^{(2^s)}(\ww)}(X)
=P_{\ww}(X)+X^{2^s}P_{\ww}(X).
\end{equation*}
\end{proposition}
%
In order to understand the nature of the evolution process in this game, let us treat the case of sequences
having a single nonzero component, which corresponds to monomials $X^k$, for some 
nonnegative integer $k$.
Let $\ww_k$ be the sequence with $\ww_k(k)=1$ and $\ww_k(j)=0$ for $j\not= k$.
Figure~\ref{Figure2} displays the first six generations produced by the game that starts with 
$\ww_3$.

By \eqref{eqpowern}, it follows that for any $m\ge 1,$
\begin{equation*}
   \begin{split}
	  P_{\psi^{(m)}}(X)&=X^k(1+X)^m\\
      &=X^k+\binom{m}{1}X^{k+1}+\cdots+\binom{m}{j}X^{k+j}+\cdots+\binom{m}{m}X^{k+m}\\
	  &=\sum_{j=0}^{m}\alpha_jX^{k+j},
   \end{split}
\end{equation*}
where $\alpha_j=1$ if $\binom{m}{j}$ is odd and $\alpha_j=0$ otherwise.
This gives a neat characterization of the entries in any generation.
%
\begin{proposition}\label{Proposition2}
   For any positive integers $k$,  $m$ and $n$,  we have
\begin{equation*}\label{eqpower2nA}
   \psi^{(m)}(\ww_k)(n)
=\begin{cases}
				  1,	\quad\text{if $\binom{m}{n-k}$ is odd}\\
				  0,	\quad\text{if $\binom{m}{n-k}$ is even}\\			
             \end{cases}
\quad\text{for any $n\in\NN$.}
\end{equation*}
\end{proposition}
%
For the numbers on the diagonal we get the following corollary.
%
\begin{corollary}\label{CorollaryTh2}
   For any $k, n\ge 1$,  we have
\begin{equation*}\label{eqpower2mn}
   \psi^{(n)}(\ww_k)(n)
=\begin{cases}
				  1,	\quad\text{if $\binom{n}{k}$ is odd}\\
				  0,	\quad\text{if $\binom{n}{k}$ is even}\\			
             \end{cases}
\quad\text{for any $n\in\NN$.}
\end{equation*}
\end{corollary}
%
The parity of the binomial coefficients leads us into the wonderful world of the Sierpinski
triangle.
Glaisher~\cite{Gla1899} proved that on any line of the Pascal arithmetic triangle, the number of
odd binomial coefficients is a power of  $2$. The reader is referred to Prunescu~\cite{Pru2012},
the authors' survey~\cite[Section 2]{CZ2013} and
the references  therein for classical and recent results on Pascal triangle modulo a power of
a prime, and their connection with Fermat numbers, constructible polygons, analysis on Sierpinski gasket,
etc.
In the context of our game, this provides us with the following result.
\begin{proposition}\label{Proposition3}
   For any integers $k\ge 1$ and $m\ge 0$, the number of nonzero entries on the line
$\psi^{(m)}(\ww_k)$
is a power of $2$.
\end{proposition}
%
\section{Action on the main side}\label{SectionMS}
A step nigher to the absolute difference game of the exponents requires reducing the evolution
field to the region from the right side of the diagonal. For this, we look for the action on the
series side.

Given a sequence  $\ww\colon\ \NN \to \{0, 1\}$, let 
$\phi(\ww)\colon\ \NN \to \{0, 1\}$ defined by 
$\phi(\ww)(n):=|\ww(n+1)-\ww(n)|$.
Then, denoting by $Q_{\ww}(X)$ the associated series 
$Q_{\ww}(X)=\sum_{n=0}^{\infty}\ww(n)X^n$, we have
\begin{equation*}
   Q_{\phi(\ww)}(X)=\frac{(1+X)Q_{\ww}(X)-Q_{\ww}(0)}{X}\,.
\end{equation*}
In order to get the subsequent generations, applying $\phi(\ww)$ repeatedly, we obtain
 \begin{equation*}
   Q_{\phi^{(m)}(\ww)}(X)=\left(\frac{1+X}{X}\right)^m Q_{\ww}(X)+ 
			   O_{\alpha\ge 1}\big(X^{-\alpha}\big)\,,
\end{equation*}
where the right term from right-hand side collects the sum monomials with negative exponents, only.
Looking for the length of the smallest possible cycle, we put $m=2^s$ and find that
 \begin{equation*}
   \begin{split}
   Q_{\phi^{(2^s)}(\ww)}(X)&=\left(\frac{1+X}{X}\right)^{2^s} Q_{\ww}(X)+ 
			   O_{\alpha\ge 1}\big(X^{-\alpha}\big)\\
      &=\frac{1+X^{2^s}}{X^{2^s}} Q_{\ww}(X)+ O_{\alpha\ge 1}\big(X^{-\alpha}\big)\\
      &=Q_{\ww}(X)+\frac{1}{X^{2^s}} Q_{\ww}(X)+ O_{\alpha\ge 1}\big(X^{-\alpha}\big)\,.%
   \end{split}
\end{equation*}
This implies that 
 \begin{equation}\label{eqQcycle}   
Q_{\phi^{(2^s)}(\ww)}(X)
      =Q_{\ww}(X)+ O_{\alpha\ge 1}\big(X^{-\alpha}\big)\,,
\end{equation}
provided that $Q_\ww(X)$ is a polynomial of degree $<2^s$.

As a consequence, relation  \eqref{eqQcycle} explains the cyclic behavior of the game generated by
the initial sequence $\ww_k$.
\begin{theorem}\label{Theorem4}
   Let $k\ge 1$ and suppose that $2^s > k$. Then, we have 
$\phi^{(2^s)}(\ww_k)=\ww_k$.
\end{theorem}
%
We remark that Theorem~\ref{Theorem4} also yields the length of the smallest cycle. 
\begin{corollary}\label{Corollary4}
   For any $k\ge 1$, the length of the smallest cycle is
 \begin{equation*}
		 L_k:=\min_{m\ge 0}\{m\colon\ m=2^s \text{ for some $s\ge 0$ and } m>k\}\,.
\end{equation*}
\end{corollary}

\noindent
{\bf Example 1.}  We provide in  Figure~\ref{Figure3} a representation of such a game
with $k=26$.  Here $s = 5$ yields the length of the smallest cycle, which is  $32$.

\section{The absolute difference game of exponents -- proof of
Theorem~\ref{Theorem1}}\label{SectionTh1}
The $Z(\cdot,\cdot)$ operation is in fact an absolute difference function on the exponents.
Working with arbitrary numbers, this is best seen by taking their localized $p$-scans.
Thus, we have  
\begin{equation*}
   Z(p^s,p^t)=p^{(s+t)-2\min(s,t)}=p^{|s-t|}.
\end{equation*}

Let $\hh : \NN\to \NN$ be the initial sequence and 
let $\eta(\hh) : \NN\to \NN$ be the transformation defined by
 \begin{equation*}
		 \eta(\hh)(n):=\frac{\hh(n)\hh(n+1)}{\big(\gcd(\hh(n), \hh(n+1)\big)^2}\,.
\end{equation*}
We shall assume that $\hh(n)$ is square free for $n\ge 0$.
Our  sequence of interest is 
$\eta^{(m)}(\hh)(0)$, $m\ge 0$.

For a given sequence $\hh$ and prime number $p$, we look at the sequence of $p$-exponents: 
$\ww_{\hh, p} : \NN\to \NN$, $\ww_{\hh, p}(n):=\nu_p(\hh(n))$ and consider the diagram

\begin{equation*}
\begin{CD}
\hh  @>     >> \ww_{\hh, p} \\
@V\eta^{(m)}VV			@VV \phi^{(m)}V\\
\eta^{(m)}(\hh)		@>  >> 	\ww_{\eta^{(m)}(\hh),p}
\end{CD}   
\end{equation*}
where 
\begin{equation}\label{eqetam}
\ww_{\eta^{(m)}(\hh),p}=\phi^{(m)}(\ww_{\hh, p}),\text{\ \ for $m\ge 1$}\,.
\end{equation}

We remark that \eqref{eqetam} is correct. Indeed, if $\nu_p(a)$ and $\nu_p(b)$ are consecutive
terms of the sequence $\ww_{\hh, p}$, then they correspond in the next generation to
$\ww_{\eta(\hh),p}$, because
\begin{equation*}
\nu_p\left(\frac{ab}{\left(\gcd(a,b)\right)^2}\right)=
\nu_p(a)+\nu_p(b)-2\nu_p\big(\gcd(a,b)\big) =
\big|\nu_p(a)-\nu_p(b) \big|\,.
\end{equation*}

Then, if $p=p_k$ is the $k$-th prime number,  by Corollary~\ref{CorollaryTh2}, it follows that 
\begin{equation*}
   \nu_{p_k}(d_m)
=\begin{cases}
				  1,	\quad\text{if $\binom{m}{k}$ is odd}\\
				  0,	\quad\text{if $\binom{m}{k}$ is even}\\			
             \end{cases}
\quad\text{for any $m\in\NN$.}
\end{equation*}
This implies that
 \begin{equation*}
   \omega(d_m)=\#\big\{k\colon\ \nu_{p_k}(d_m)=1\big\}\\
		=\#\big\{k\colon\ \binom{m}{k} \text{ is odd}\big\}\,,	
\end{equation*}
which  is a power of $2$, by a result of Glaisher~\cite{Gla1899}. This proves
Theorem~\ref{Theorem1}.

\section{The game revisited -- proofs of Theorems~\ref{Theorem1}
and~\ref{Theorem2}}\label{SectionTh1Th2}

In previous sections we saw that one gets the best view of the field of the game on exponents
by considering it as a superposition of $p$-renderings for all prime numbers $p$.
In Figure~\ref{Figure3}, there is an example of such a rendering for $p=101$, the $26$-th prime. In
Section~\ref{SectionMS} we have noticed that incomplete Sierpinski triangles repeat in a cycle of
length $32$ (the first number that is a power of $2$ and is larger than $26$). 
This is just an instance of an arbitrary prime and there is a rewarding way to look on it. Indeed,
let $T_{\infty}$ be an infinite Sierpinski triangle with one vertex on the top
and having cells corresponding to the honeycomb network of our game. 
We assume that $T_{\infty}$ runs across the honeycomb moving from left to right, one step at a time.
Furthermore, we assume that triangle $T_{\infty}$  changes to a new color its nonempty
cells at every step.
\bigskip

\begin{center}
\begin{figure}[ht]
\centering
   \includegraphics[scale=0.93,angle=-90]{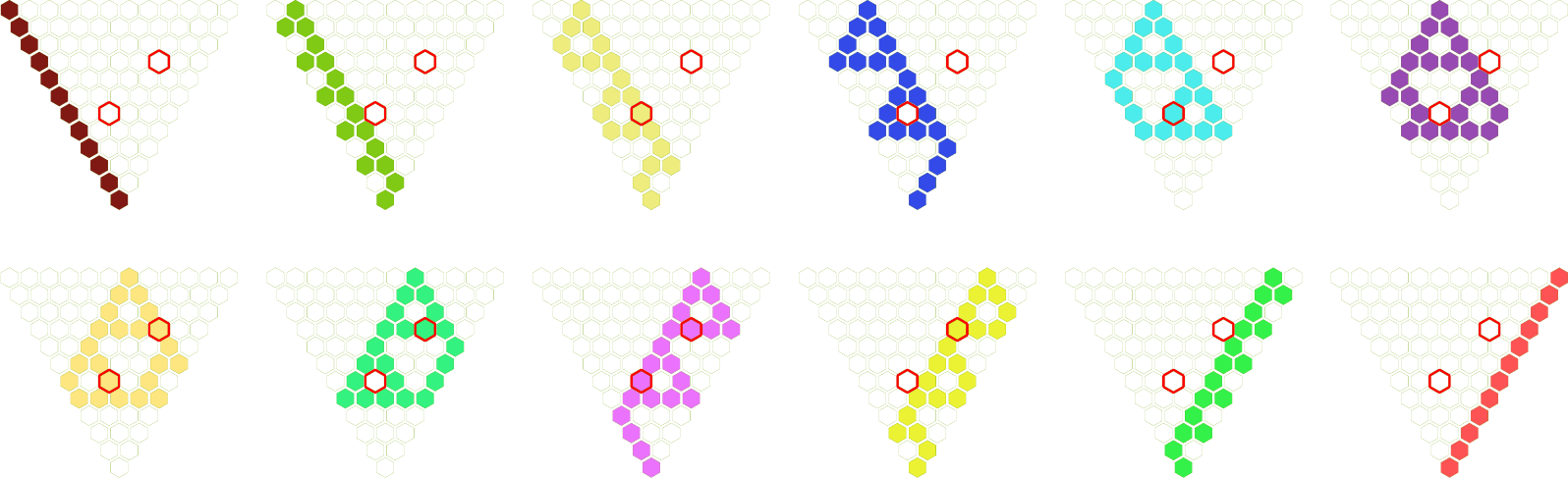}
\bigskip

\caption{Cut-off triangles of size $12$ from the the top-left of the field game covered by
the pictures
$\PP_1$, $\PP_2$, \ldots, $\PP_{12}$.
On each of them, focusing the sight, are highlighted the cells that correspond to the numbers
$a_{3, 7}$ and $a_{6, 3}$.}
 \label{Figure1-12}
 \end{figure}
\end{center}

\begin{center}
\begin{figure}[ht]
\centering
   \includegraphics[width=0.67\textwidth, angle=-90]{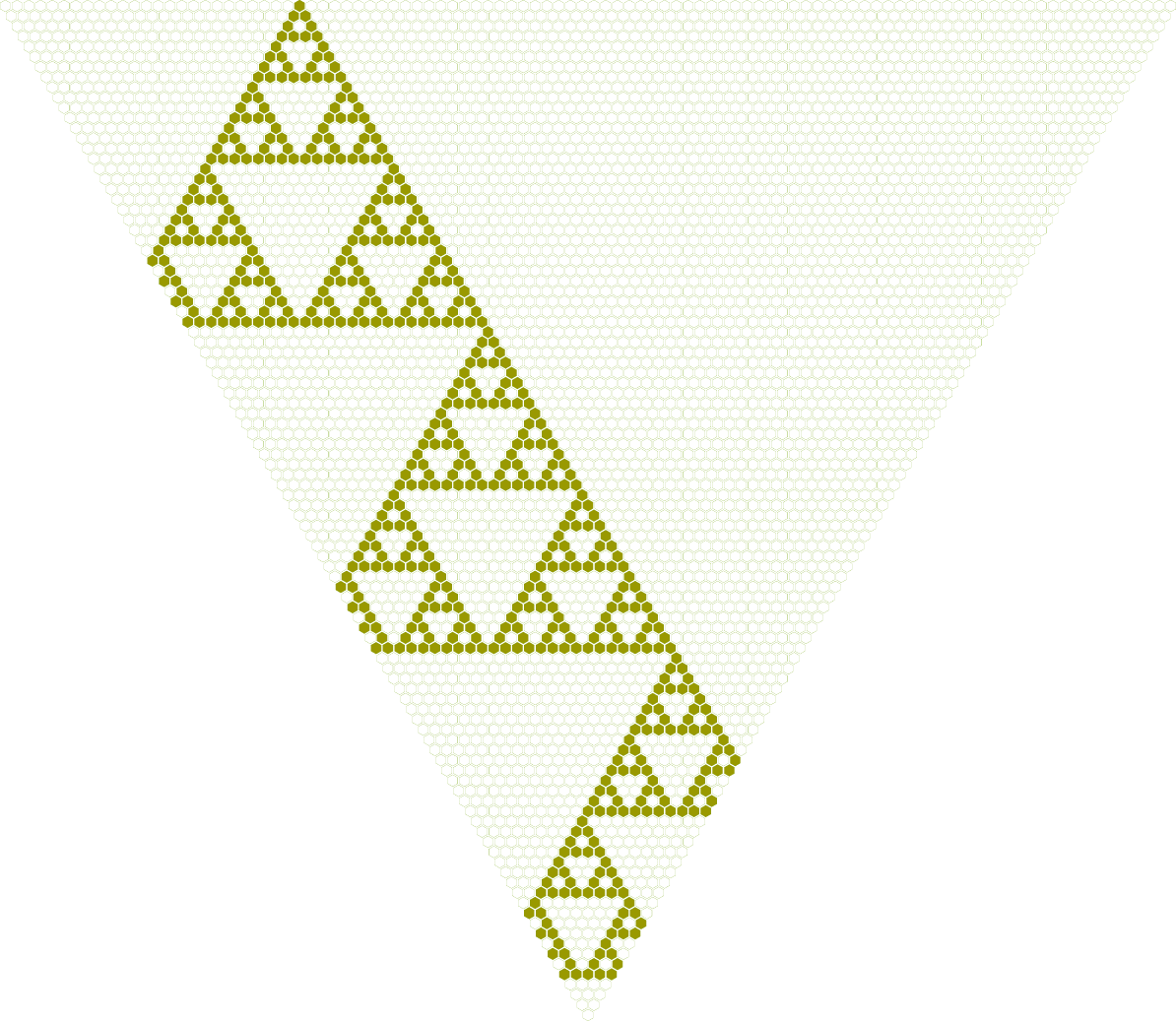}
\caption{Triangular slice of size $N=100$ cut off from
the rendering produced by $p=101$, the  $26$-th prime, in the table generated by  iterated
application of $Z(\cdot, \cdot)$ on the sequence of prime numbers.
}
 \label{Figure3}
 \end{figure}
\end{center}

For $n\ge 1$, we denote by $\PP_n$ the image of $T_{\infty}$ that lies over the honeycomb of the
game when its vertex from the top lies over the $n$-th cell of the
topmost line.
Thus, $\PP_n$ is the $p$-rendering of the game corresponding to $p=p_n$, the $n$-th prime.
For example, in Figures~\ref{Figure1-12} and~\ref{Figure3}
 are shown the upper parts of $\PP_1$, $\PP_2$, \ldots,
$\PP_{12}$ and $\PP_{26}$, respectively.
With these realizations, our number game from Figure~\ref{FigureZprimes10} corresponds to the
mosaic of colors developed by the superposition of all $\PP_n$, $n\ge 1$.

Now let us look at an arbitrary cell   $m\ge 0$,  $n\ge 1$. 
The geometrical constrains show that the prime numbers $p_k$, $k\ge 1$, for which a nonempty cell of
$\PP_k$ lies over the cell of $a_{m,n}$ must satisfy the inequality 
$p_n\le p_k\le p_{m+n}$. 
Moreover, we need to consider only those primes $p_k$ for which a nonempty cell of $\PP_k$
lies over the cell $(m,n)$. 
Let $\mathcal{S}_m$ be the set of nonzero numbers from the $m$-th row of Sierpinski triangle, that
is,
\begin{equation*}
   \mathcal{S}_m:=\Big\{r\colon\ \binom{m}{r}\equiv 1\pmod 2, \ \ 0\le r\le m\Big\}\,.
\end{equation*}
Then, we have obtained a closed form expression for any entry in Figure~\ref{FigureZprimes10}. This
is
\begin{equation*}
   a_{m,n}=\prod_{r\in\mathcal{S}_m}p_{n+r}\,.
\end{equation*}
In particular, we have proved both Theorem~\ref{Theorem1} and
 Theorem~\ref{Theorem2}. 

\section{Conclusion and Future Directions}
We discussed the outcome of the iterated application of a simple arithmetical function $Z(a,b)$,
which operates like a recorder of the result produced by the  plastic collision of two elementary
particles. Thus $Z(a,b)$ is a product of powers of primes $p^{d_p}$ and locally, for any prime $p$, 
the exponent $d_p$ is the absolute difference of the exponents of $p$ in $a$ and $b$.
This paper leads to some further directions of research.
Some open problems that naturally arise are presented below.

Sequence $d_m$ from Table~\ref{Table1} increases in dyadically enhanced twists. 
It combines the size of primes with the selection process structured by the Sierpinski triangle
sieve.
Therefore, the following natural questions arise.

\bigskip

\noindent
\textbf{Question 1.} Find accurate lower and upper bounds for

\begin{equation*}
  	 \min_{x\le m\le y} d_m,\ \ 
	 \min_{x\le m\le y}\omega(d_m),\ \ 
	 \max_{x\le m\le y} d_m\ \  \text{and}\ \  
   \max_{x\le m\le y}\omega(d_m). 
\end{equation*}

\smallskip

\noindent
\textbf{Question 2.} Estimate the following averages:
 \begin{equation*}
      \sum_{x\le m\le y}\omega(d_m) \ \ \text{ and }\ \ 
	  \sum_{x\le m\le y}d_m.
\end{equation*}

\noindent
\textbf{Question 3.}
Starting with the sequence of positive integers instead of the primes on the first line and applying
the same generating 
$Z(\cdot, \cdot)$- rule, we obtain a triangle that is similar to the triangle in
Figure~\ref{FigureZprimes10} (see\cite[Section 9]{CZ2013}).
%
The sequence of numbers on the left side of this new triangle is 
(see  \cite[sequence {\tt A222311}]{OEIS})
\begin{equation}\label{eqsnat}
 1, 2, 3, 6, 5, 15, 105, 70, 1, 5, 33, 55, 65, 273, 1001,\ldots
\end{equation}
By sorting it and removing duplicates we get a new sequence. 
Sloane \cite[sequence {\tt A222313}]{OEIS} makes this operations with  the first $500$ terms of
\eqref{eqsnat} to obtain:
$  1, 2, 3, 5, 6, 15, 17$, $33, 55, 65, 70, 105, \ldots$
The first three terms shown are certain. 
Sloane asks how far this list is complete up to $105$ and if there is a proof that $4$ cannot
appear.


\smallskip\smallskip\smallskip\smallskip
\noindent
\textbf{Acknowledgement}
\smallskip\smallskip

\noindent
The authors would like to thank the anonymous referees for their valuable comments and
suggestions, which helped improve the presentation of the paper.  


\smallskip\smallskip\smallskip\smallskip\smallskip\smallskip\smallskip\smallskip

%
%

\end{document}